\documentclass{amsart}
\usepackage{amsmath,amsthm,amsfonts,amssymb, verbatim}
\usepackage{graphicx,subfigure}
\usepackage{hyperref}
\usepackage[initials,nobysame]{amsrefs}

\iftrue
\newcommand{\sidenote}[1]{\marginpar[\raggedleft\tiny #1]{\raggedright\tiny #1}}
\else
\newcommand{\sidenote}[1]{}
\fi

\newcommand{\del}{\partial}

\newcommand{\lap}{\triangle}

\newcommand{\inv}{^{-1}}

\newcommand{\grad}{\nabla}

\newcommand{\divergence}{\grad \cdot}

\newcommand{\commentout}[1]{}
\newcommand{\eps}{\epsilon}
\newcommand{\pdr}[2]{\frac{\partial{#1}}{\partial{#2}}}

\renewcommand{\epsilon}{\varepsilon}
\renewcommand{\leq}{\leqslant}
\renewcommand{\geq}{\geqslant}

\newcommand{\gradperp}{\grad^\perp}

\newcommand{\R}{\mathbb{R}}

\newif\iftextstyle
\textstyletrue
\everydisplay\expandafter{\the\everydisplay\textstylefalse}
\newcommand{\ds}{\displaystyle}

\newcommand{\abs}[1]{\iftextstyle\lvert#1\rvert\else\left\lvert#1\right\rvert\fi}
\newcommand{\norm}[1]{\iftextstyle\lVert#1\rVert\else\left\lVert#1\right\rVert\fi}

\newcommand{\leb}[1]{L^#1}
\newcommand{\lpnorm}[2]{\norm{#1}_{\leb{#2}}}
\newcommand{\linfnorm}[1]{\lpnorm{#1}{\infty}}

\newcommand{\at}{\iftextstyle|\else\Big|\fi}

\numberwithin{equation}{section}
\allowdisplaybreaks

\newtheorem{theorem}{Theorem}[section]
\newtheorem{lemma}[theorem]{Lemma}
\newtheorem{proposition}[theorem]{Proposition}

\newtheorem*{theorem*}{Theorem}
\newtheorem*{lemma*}{Lemma}
\newtheorem*{proposition*}{Proposition}
\newtheorem*{corollary*}{Corollary}

\theoremstyle{definition}

\theoremstyle{remark}
\newtheorem{remark}[theorem]{Remark}
\newtheorem*{remark*}{Remark}

\begin{document}
\title{Exit times of diffusions with incompressible drift}
\author{Gautam Iyer}
\address{Department of Mathematical Sciences, Carnegie Mellon University, Pittsburgh PA 15213}
\email{gautam@math.cmu.edu}
\author{Alexei Novikov}
\address{Department of Mathematics, Pennsylvania State University, State College PA 16802}
\email{anovikov@math.psu.edu}
\author{Lenya Ryzhik}
\address{Department of Mathematics, Stanford University, Stanford CA 94305}
\email{ryzhik@math.stanford.edu}
\author{Andrej Zlato\v s}
\address{Department of Mathematics, University of Chicago, Chicago IL 60637}
\email{zlatos@math.uchicago.edu}
\subjclass[2000]{%
35J60, 
35J05. 
}
\begin{abstract}
Let  $\Omega\subset\R^n$ be a bounded domain and for $x\in\Omega$ let $\tau(x)$ be the expected exit time from  $\Omega$ of a diffusing particle starting at $x$ and advected by an incompressible flow $u$. We are interested in the question which flows maximize $\norm{\tau}_{L^\infty(\Omega)}$, that is, they are most efficient in the creation of hotspots inside $\Omega$. Surprisingly, among all simply connected domains in two dimensions, the discs are the only ones for which the zero flow $u\equiv 0$ maximises $\norm{\tau}_{L^\infty(\Omega)}$. We also show that in any dimension, among all domains with a fixed volume and all incompressible flows on them, $\norm{\tau}_{L^\infty(\Omega)}$ is maximized by the zero flow on the ball.
\end{abstract}
\maketitle
\section{Introduction}
 
It is well-known that mixing by an incompressible flow enhances   diffusion in many contexts. This is demonstrated, for instance, by the fact that the effective diffusivity
of a periodic incompressible flow is always larger than diffusion in the absence of a flow~\cite{FannP}, or that the principal
eigenvalue $\mu_u$ of the problem
\begin{equation}\label{aug251}
\left\{
\begin{array}{cl}
-\lap\phi+u\cdot\nabla\phi=\mu_u\phi & \text{ $\phi>0$ in $\Omega$},\\
\phi=0 & \text{ on $\partial\Omega$}
\end{array}\right.
\end{equation}
is never smaller than the corresponding eigenvalue $\mu_0$ of \eqref{aug251} with $u\equiv 0$. Classes of flows which are most effective in enhancing diffusion have been studied both on bounded and unbounded domains, and their characterizations have been provided in~\cites{CKRZ, Zla2D}.

On the other hand, it was observed in~\cite{BKJS} that an incompressible flow may actually 
slow down diffusion in the following sense. Consider the explosion problem
$$
\begin{array}{cl}
-\lap\phi+u\cdot\nabla\phi=\lambda e^{\phi} & \text{ in $\Omega$},\\
\phi=0 & \text{ on $\partial\Omega$.}
\end{array}
$$
There exists $\lambda_*(u)$ such that this problem has a solution for all $\lambda\leq\lambda_*(u)$ and no solution for $\lambda>\lambda_*(u)$ (see  \cites{CR,JL,KK} for $u\equiv 0$ and~\cite{BKNR} for $u\not\equiv 0$). Surprisingly, it was shown {\it numerically} in~\cite{BKJS} that in a long rectangle there are incompressible flows with $\lambda_*(u)<\lambda_*(0)$. This means that addition of a flow (which typically increases $\lambda_*$ due to mixing) can sometimes instead promote the creation of hotspots and inhibit their interaction with the cold boundary $\partial\Omega$.

The present paper is a step toward mathematical understanding of this diffusion slowdown effect of certain incompressible flows. We consider the  problem
\begin{equation}\label{eqnEExitTime}
\left\{
\begin{array}{cl}
- \lap \tau^u+\ds u \cdot \grad \tau^u   = 1 &\text{in $\Omega$,}\\
\ds \tau^u = 0 & \text{on $\del \Omega$}
\end{array}
\right.
\end{equation}
on a smooth bounded domain $\Omega \subset \R^n$, with $u(x)$ an incompressible flow on $\Omega$ (i.e., $\nabla\cdot u\equiv 0$) which is tangential to $\partial\Omega$ (i.e., $u\cdot\hat n\equiv 0$ on $\partial\Omega$, with $\hat n$ the outward normal to $\partial\Omega$).
Physically, the solution $\tau^u(x)$ is the expected exit time from $\Omega$ of the random process 
\[
dX_t=-u(X_t)\,dt+\sqrt{2}\,dB_t,\qquad X_0=x,
\]
modeling the motion of a diffusing particle advected by the flow $u$. Although one might think that the expected exit time is always decreased by the addition of an incompressible flow due to improved mixing, this need not be the case. Our first result shows that in any bounded simply connected domain in $\R^2$ which is not a disk, there are (regular) incompressible flows which \textit{increase} the maximum of the expected exit time of $X_t$ from $\Omega$.

\begin{theorem}\label{thmWorseOnNonDisk}
Let $\Omega \subset \R^2$ be a bounded simply connected domain with a $C^1$ boundary which is not a disk. Then there exists a $C^1$, divergence free vector field $v:\Omega \to \R^2$  tangential to $\del\Omega$ such that $\norm{\tau^v}_{L^\infty(\Omega)} > \norm{\tau^0}_{L^\infty(\Omega)}$.
\end{theorem}

{\it Remark.}
We note that the incompressible flows are the natural class to study in this context. Indeed, if one considers general $v$ (not necessarily divergence free), then it is  easy to show that $\norm{\tau^v}_{L^\infty(\Omega)}$ can be made arbitrarily large by, for instance, taking $v(x)=A(x_0-x)$ with $x_0\in\Omega$ and $A$ sufficiently large.
\smallskip

On a disk however, no (incompressible) stirring will increase this expected exit time beyond the one for $u\equiv 0$. In fact we prove in any dimension that the $L^p$-norm of the expected exit time can never be larger than that from a disk of equal volume with $u\equiv 0$.

\begin{theorem}\label{T.1.1}
Let $\Omega\subset \mathbb{R}^n$ be a bounded domain with a $C^1$ boundary and $v:\Omega \to \R^n$ a $C^1$ divergence free vector field tangential to $\del \Omega$. 
Then for any $p\in[1,\infty]$,
$$
\norm{\tau^v}_{L^p(\Omega)} \leq \norm{\tau^{0,D}}_{L^p(D)}
$$
where $D \subset \R^n$ is a ball with the same Lebesgue measure as $\Omega$, and $\tau^{0,D}$ is the solution of~\eqref{eqnPoisson} on $D$ with $u \equiv 0$.
\end{theorem}

{\it Remark.}
If $D$ is a ball with Lebesgue measure $V$ and center $0$, then $\tau^{0,D}$ is given explicitly by the formula
$$
\tau^{0,D}(x) = \frac{1}{2n}\left[ \left(\frac{V}{\Gamma_n}\right)^\frac{2}{n}-\abs{x}^2 \right],
$$
with $\Gamma_n$ the Lebesgue measure of the unit ball in $\mathbb{R}^n$.
\smallskip

There are, of course, other ways to quantify the effect of stirring on diffusion --- see, for instance, \cite{STD} where many additional references can be found, especially to the physics literature. Closely related to the problem studied in the present paper is the following question. It is shown in~\cite{BKNR} that for any $p>n/2$ there exists
a constant $C_p(\Omega)$ such that for any incompressible $u$ tangential to $\partial\Omega$ and any $f\in L^p(\Omega)$, the solution of
\begin{equation}\label{eqnPoisson}
\left\{\begin{array}{cl}
- \lap \phi+\ds u \cdot \grad\phi  = f &\text{in $\Omega$,}\\
\ds \phi = 0 & \text{on $\del \Omega$}
\end{array}\right.
\end{equation}
satisfies $\norm{\phi}_{L^\infty}\leq C_p(\Omega)\norm{f}_{L^p}$. It would be interesting to determine which flows achieve $C_p(\Omega)$ and how does $C_p(\Omega)$ depend on $\Omega$. Theorems~\ref{thmWorseOnNonDisk} and \ref{T.1.1} are a first step in this direction.

The present paper is organized as follows. In Section~\ref{proofs} prove our main results, Theorems~\ref{thmWorseOnNonDisk} and~\ref{T.1.1}. Our proof of Theorem~\ref{thmWorseOnNonDisk} involves a variational principle, Proposition~\ref{ppnCriticalPt}, the proof of which is somewhat technical and therefore postponed to Section~\ref{ptproof}. This variational principle leads to an interesting PDE for the critical points of the expected exit time functional. We discuss properties of these critical points and provide some numerical examples in Section~\ref{maximizer}.

\subsection*{Acknowledgment}
This work was supported in part by NSF grants of the authors. 
AZ was also supported by an Alfred P.~Sloan Research Fellowship.  
GI was also partially supported by the Center for Nonlinear Analysis.

\section{Proofs of the main results}\label{proofs}
\subsection{Proof of Theorem~\ref{thmWorseOnNonDisk} }
For a given incompressible $C^3$ flow $u$ tangential to $\partial\Omega$, 
consider the family of Poisson problems
\begin{equation}\label{aug181}
\left\{\begin{array}{cl}
-\lap\tau^{Au}+Au\cdot\nabla\tau^{Au}=1 & \text{ in $\Omega$},\\
\tau^{Au}=0 & \text{ on $\partial\Omega$}
\end{array}\right.
\end{equation}
with $A>0$. Let $\psi$ be the {\it stream function} of $u$, that is, $\psi:\Omega\to\R$ is $C^4$ and such that $\psi(\partial\Omega)=0$ and  $u=\nabla^\perp\psi=(-\del_{2} \psi, \del_{1} \psi)$.
It is well-known~\cites{Freidlin,WF} that if all critical points of  $\psi$ are
non-degenerate and no two of them lie on the same level set of $\psi$, then the functions $\tau^{Au}$ converge uniformly to a limit $\bar\tau^u$ which is constant on the level sets of $\psi$ and satisfies an asymptotic Freidlin problem on the Reeb graph of the function $\psi$. If $\Omega$ is simply connected and $\psi$ has a single (non-degenerate) critical point 
(in which case either $\psi>0$ or $\psi<0$ on $\Omega$, and we will assume without loss the former), then we have the explicit formula
\begin{equation}\label{aug182}
\bar \tau^u(y) = - \int_0^{\psi(y)} \frac {\abs{\Omega_{\psi,h}} } {\int_{\Omega_{\psi,h}} \lap\psi \, dx} \, dh.
\end{equation}
Here and elsewhere we let $\Omega_{\psi,h}= \{x\in\Omega\,\big|\,\psi(x)>h\}$, the $h$-super-level set of $\psi$. Notice that $\bar \tau^u$ is just a reparametrization of $\psi$.
As we prove in Proposition~\ref{ppnFreidlinProblem} below, the formula \eqref{aug182} holds also when the single critical point of $\psi$ is degenerate.

We will start by considering only flows with the above property. That is, $u$ is $C^3$ and such that the stream function $\psi$ only has a single critical point in $\Omega$ (which is simply connected
and $\psi>0$ on $\Omega$). In particular, all super-level sets of $\psi$ are simply connected and $\psi$ attains a single maximum $\psi(x_0)=M>0$. Moreover, any $h\in(0,M)$ is a regular value of $\psi$ and $\partial\Omega_{\psi,h}$ is a $C^4$ Jordan curve. 

Assume now that for some $C^3$ incompressible flow $w$ the function $\tau^w$ has a single critical point and let $\psi=\tau^w$ (which is $C^4$) and $u=\grad^\perp\psi$. Thus $\psi$ solves 
\begin{equation}\label{aug183}
\left\{\begin{array}{cl}
-\lap\psi+w\cdot\nabla\psi=1 & \text{ in $\Omega$},\\
\psi=0 & \text{on $\partial\Omega$}.
\end{array}\right.
\end{equation}  
Integrating this over $\Omega_{\psi,h}$ and using incompressibility of $w$, we obtain 
\[
\abs{\Omega_{\psi,h}}= -\int_{\Omega_{\psi,h}} \lap\psi \, dx
\]
 for any $h\in(0,M)$. This  together with (\ref{aug182})
 implies that $\bar\tau^u \equiv\psi$. That is, such solutions $\psi$ to the Poisson problem \eqref{aug183} solve the Freidlin problem for themselves. We are particularly interested in the case $w=0$, with $\psi=\tau^0$ solving
\begin{equation}\label{aug184}
\left\{\begin{array}{cl}
-\lap\tau^0=1 & \text{ in $\Omega$},\\
\tau^0=0 & \text{ on $\partial\Omega$}
\end{array}\right.
\end{equation}
and  $u_0=\nabla^\perp\tau^0$. Notice that then $\tau^0$ also solves
\begin{equation}\label{aug186}
\left\{\begin{array}{cl}
-\lap\tau^0+Au_0\cdot\nabla\tau^0=1 & \text{ in $\Omega$},\\
\tau^0=0 & \text{ on $\partial\Omega$.}
\end{array}\right.
\end{equation}
for any $A\in\R$ and so $\tau^0=\bar\tau^{u_0}$. 
Let us therefore assume, for now, that $\Omega$ is such that $\tau^0$ has a single 
critical point. 

We now assume that for any incompressible flow $u$ on $\Omega$ we have $\norm{\tau^u}_{L^\infty}\le \norm{\tau^0}_{L^\infty}$. In particular, $\norm{\bar\tau^u}_{L^\infty}\le \norm{\bar\tau^{u_0}}_{L^\infty}$ for each $u$ whose stream function has a single critical point. We will now show that this is the case only when $\Omega$ is a disc, thus proving Theorem~\ref{thmWorseOnNonDisk} for all $\Omega$ such that $\tau^0$ has a single critical point. 

The key ingredient of our proof is that for all ``infinite amplitude'' expected exit times $\bar\tau^u$, we have a variational principle which gives an explicit equation satisfied by the critical points (and thus the maximiser) of the functional $I(\psi)=\|\bar\tau^{\grad^\perp\psi}\|_{L^\infty}$ (with $\psi$ having a single critical point).

We set up the variational principle as follows. Let $v$ (the ``direction'' of our variation) be any 
$C^4$ vector field tangential to $\del \Omega$. Let $X$ be the flow (in the dynamical systems sense) given by
$$
\frac{dX_\epsilon}{d\epsilon} = v\circ X_\epsilon, \qquad X_0 = \text{Id}.
$$
Given a stream function $\psi$ with a single critical point, we perturb it by composing it with the flow $X_\epsilon$. 
Let $\psi^\epsilon = \psi \circ X_\epsilon$, $u^\epsilon = \gradperp \psi^\epsilon$, and $\bar \tau^\epsilon = \bar \tau^{u^\epsilon}$. 
Notice that $\psi^\epsilon$ is $C^4$ and again has a single critical point (the maximum) $x_0^\epsilon$. Then $\bar\tau^\eps$ also attains its maximum at $x_0^\epsilon$ due to \eqref{aug182}, so the variation of $I$ in direction $v$ is
\begin{equation}\label{eqnVariationTau}
V(\psi, v ) =
\left.\frac{d}{d\epsilon} \bar \tau^\epsilon(x_0^\epsilon) \right\rvert_{\epsilon = 0}.
\end{equation}
We say that  $\psi$ is a critical point of $I$ if for all $C^4$ (not necessarily divergence free) vector fields $v$ tangential to $\del \Omega$, we have  $V(\psi, v )=0$. Clearly any  $\psi$ (with a single critical point) which maximises  $I$  is a critical point of $I$. So our aim is to prove that $\tau^0$ is not a critical point of $I$ unless $\Omega$ is a disc (assuming for now that $\tau^0$ has a single critical point).

As mentioned earlier, the proof of this fact rests on obtaining an explicit equation for critical points of $I$. 
We can now do this   by a direct computation using the Freidlin-Wentzel theory~\cites{WF,Freidlin}.

\begin{proposition}\label{ppnCriticalPt}
Let $\Omega\subset\R^2$ be a bounded simply connected domain with a $C^1$ boundary and let $\psi>0$ be a $C^4$ stream function on $\Omega$ with $\psi(\partial\Omega)=0$ and a single 
critical point. Then $\psi$ is a critical point of the functional $I$ if and only if $\phi = \bar \tau^{\gradperp \psi}$, the solution of the Freidlin problem~\eqref{aug182} with stream function $\psi$, also
solves
\begin{equation}\label{eqnCriticalPt}
-2\lap \phi(x) = 1 +  {\abs{\grad \phi(x)}^2} 
\int_{\partial\Omega_{\phi,\phi(x)}} \frac{d\sigma}{\abs{\grad \phi}} 
{\ds\left(\int_{\partial\Omega_{\phi,\phi(x)}} \abs{\grad \phi} \, d\sigma \right)^{-1}} 
\end{equation}
\end{proposition}

We postpone the proof of Proposition~\ref{ppnCriticalPt} to Section~\ref{ptproof}, but make two remarks before proceeding.

\begin{remark*}
One can also write down an explicit PDE \eqref{eqnCriticalPtVarphi} for the stream function $\psi$. This PDE, however, is somewhat more complicated, and we find it more convenient to work with~\eqref{eqnCriticalPt} involving the reparametrization $\phi$ of $\psi$.
\end{remark*}

\begin{remark*}
Assume that a $C^3$ flow $w$ maximizes $\|\tau^w\|_{L^\infty}$ and $\tau^w$ has a single critical point. The argument following \eqref{aug183} above then shows  $\tau^w\equiv \bar\tau^{\nabla^\perp\tau^w}$, so $\tau^w$ is a critical point of $I$ and solves \eqref{eqnCriticalPt}.
\end{remark*}


Let now $\psi=\tau^0$ have a single critical point $x_0\in\Omega$ and assume that $\psi$ is a critical point of $I$. Recall that $\psi=\bar\tau^{u_0}=\bar\tau^{\nabla^\perp \psi}$, so  Proposition~\ref{ppnCriticalPt} implies that $\psi$ solves~\eqref{eqnCriticalPt}.  Since $-\lap \psi = 1$, we obtain
\begin{equation}\label{eqnReducesToEikonal}
\abs{\grad \psi(x)}^2 \int_{\partial\Omega_{\psi,\psi(x)}} \frac{d\sigma}{\abs{\grad \psi}}\left(\ds\int_{\partial\Omega_{\psi,\psi(x)}} \abs{\grad \psi} \, d\sigma \right)^{-1} =1,
\end{equation}
immediately showing that $\abs{\grad \psi}$ must be constant on the level sets of $\psi$. 
Thus $\psi$ solves the eikonal equation $\abs{\grad \psi(x)}=g(\psi)$ with $g$ equal zero at the maximum of $\psi$ and positive elsewhere. 
It is well known that a solution of such 
equation does not have interior singularities only if  $\Omega$ is a disk and $\psi$ 
is radial~\cite{BBI}. 
In our situation this can be seen as follows. After reparametrization we may assume that $g\equiv1$, and $\psi$ attains its maximum at $x_0$. This introduces a singularity at $x_0$ so let us suppose that $\psi$ does not have other interior singularities. Since the level sets of $\psi$ are connected, and the maximum is isolated, for any $\epsilon>0$ we can find a wavefront (a level set of $\psi$) that is contained in a disc of radius  $\epsilon>0$ centered at $x_0$. By compactness, this wavefront is a positive
distance $\epsilon' \in (0, \epsilon)$  away from $x_0$. Absence of singularities now implies that we can evolve this level set, and the spheres of radius $\epsilon$ and $\epsilon'$ ``outward'' by the eikonal equation (with $g\equiv 1$). Then each level set of $\psi$ obtained by this evolution lies entirely within distance $\epsilon - \epsilon' < \epsilon$ from a circle. As $\epsilon>0$ is arbitrary we conclude that level sets of $\psi$ have to be circles. Since $ \psi(\partial\Omega) =0$, we have that $\Omega$ is a disk and $\psi$ radial.

Thus we have proved that if $\tau^0$ has a single critical point, then it does not maximize $I$ when $\Omega$ is not a disc. Since the claim of Theorem~\ref{thmWorseOnNonDisk} for a disc follows from Theorem~\ref{T.1.1} (which we will prove shortly), we are left with considering the case of $\Omega$ such that $\tau^0$ has more than one critical point. We will use the following claim to reduce this to the previous case.

\begin{lemma}
For any bounded simply connected $\Omega\subset\R^2$ with a $C^1$ boundary,  the set of maxima of $\tau^0$ is discrete. 
\end{lemma}

\begin{proof}
Let $M = \| \tau^0\|_{L^\infty}$, let $\mathcal D = \{x\,|\,\tau^0(x)=M\}$, and suppose $\mathcal D$ is not discrete. Since $\mathcal D$ is positive distance from $\partial\Omega$, it has an accumulation point $x_0$ inside $\Omega$. Assume without loss of generality
that a sequence $x_n\in\mathcal D$ converges to $x_0$ along the $x$-axis: $(x_n-x_0)/|x_n-x_0|\to (1,0)$. Thus $\partial_x\tau^0(x_0)=\partial_x^2\tau^0(x_0)=0$, so $\partial_y^2\tau^0(x_0)=1$ and the analytic implicit function theorem shows that there is a real analytic curve $\mathcal C$ containing $x_0$ on which $\partial_y\tau^0=0$ (since $\tau^0$ is real analytic). It then follows that $x_n\in\mathcal C$ for all large $n$, and real analyticity of $\tau^0|_{\mathcal C}$ now shows $\mathcal C\subset \mathcal D$.

So $\mathcal D$ contains an analytic curve which cannot end inside $\Omega$ (by the previous argument) and must also stay away from $\del \Omega$ (by $M>0$). This means that such a curve must be closed. But then $\tau^0 > M$ inside the region enclosed by this curve (which is a subset of $\Omega$), a contradiction.
\end{proof}

We now return to the proof of Theorem~\ref{thmWorseOnNonDisk} for general $\tau^0$ and assume that the zero flow maximizes $\|\tau^u\|_{L^\infty}$ among all incompressible flows $u$ on $\Omega$. We will reduce the problem to the previous case by showing that then the same is true for a connected component of $\Omega_{\tau^0,h}$ containing a maximum of $\tau^0$, for all $h$ sufficiently close to $M = \linfnorm{\tau^0}$. We introduce some notation and make this precise below.

Let $\tau^0(x_0)=M$ for some $x_0\in\Omega$ and denote by
$\Omega_h$ the connected component of $\Omega_{\tau^0,h}$ containing $x_0$.  For any $h$ and incompressible $C^3$ vector field $w$ tangential to $\del\Omega_h$, define $Q_{\Omega_h}(w) = \limsup_{A\to\infty} \norm{\tau^{Aw}_{\Omega_h}}_{L^\infty(\Omega_h)}$ (where $\tau^{Aw}_{\Omega_h}$ satisfies~\eqref{eqnEExitTime} with $\Omega = \Omega_h$ and $u = Aw$). Finally,  choose $h_0<M$ sufficiently close to $M$ so that $\bar\Omega_{h_0}$ contains no critical points of $\tau^0$ besides $x_0$.

\begin{lemma}\label{lemaug212}
Assume that for all $C^3$ incompressible vector fields $u$ tangential to $\del \Omega$, we have $\norm{\tau^0}_{L^\infty(\Omega)} \geq Q_{\Omega}(u)$. Then for any  $C^3$ incompressible vector field $w$  tangential to $\del \Omega_{h_0}$, we have $\norm{\tau^0-h_0}_{L^\infty(\Omega_{h_0})} \geq Q_{\Omega_{h_0}}(w)$.
\end{lemma}

Momentarily postponing the proof of Lemma~\ref{lemaug212}, note that $\tau^0 - h_0$ is the expected exit time of Brownian motion, starting at $x$, from $\Omega_{h_0}$. That is, $\tau^0 - h_0$ is the solution of~\eqref{eqnEExitTime} with $\Omega = \Omega_{h_0}$, and $u = 0$. Thus Theorem~\ref{thmWorseOnNonDisk} for $\Omega_{h_0}$ (which we have already proved) shows that $\Omega_{h_0}$ is a disk of some radius, say $R$, and $\tau^0$ is radial in it. Since $\tau^0(x)+\tfrac 12|x-x_0|^2$ is harmonic in $\Omega$ and radial near $x_0$, it must be constant in $\Omega$. Since $\tau^0(\partial\Omega)=0$ we have that $\Omega$ is a disk, completing the proof of Theorem~\ref{thmWorseOnNonDisk}.

It only remains to prove Proposition~\ref{ppnCriticalPt} and Lemma~\ref{lemaug212}. Proposition~\ref{ppnCriticalPt} is proved in Section~\ref{ptproof}, and we prove Lemma~\ref{lemaug212} below.

\begin{proof}[Proof of Lemma~\ref{lemaug212}]
The proof is based on the more general observation that changing any stream function near its maximum does not affect the asymptotic $A\to\infty$ behavior of the solution of~\eqref{aug181} away from the maximum. We make this precise below.
 
Let $\psi$ be any $C^4$ function in $\Omega$ with $\psi(\partial\Omega)=0$ and let $x_0,M,h_0,\Omega_h$ be defined as above, with $\psi$ in place of $\tau^0$ (we will eventually choose $\psi = \tau^0$).
For some $h_1\in(h_0,M)$ let $\psi'$ be some $C^4$  function such that $\psi(x)=\psi'(x)$ for $x\in\Omega\setminus\Omega_{h_1}$, and denote $u=\nabla^\perp \psi$, $u'=\nabla^\perp \psi'$. Let $\tau_A$ and $\tau_A'$ solve
\begin{equation}\label{aug213}
\left\{\begin{array}{cl}
-\lap\tau_A+Au\cdot\nabla\tau_A=1 & \text{ in $\Omega$},\\
\tau_A=0 & \hbox{ on $\partial\Omega$,}
\end{array}\right.
\end{equation}
and
\begin{equation}\label{aug214}
\left\{\begin{array}{cl}
-\lap\tau_A'+Au'\cdot\nabla\tau_A'=1 &\text{ in $\Omega$},\\
\tau_A'=0 &\text{ on $\partial\Omega$.}
\end{array}\right.
\end{equation}
We will first show
\begin{equation}\label{aug216}
\norm{\nabla\tau_A-\nabla\tau_A'}_{L^2(\Omega\setminus\Omega_{h_0})}\to0\text{ as }A\to+\infty,
\end{equation}
which, as mentioned earlier, says that perturbations of the stream function near $x_0$ do not affect the asymptotic $A\to\infty$ behavior away from $x_0$.

To prove~\eqref{aug216}, let $\phi_A=\tau_A-\tau_A'$ so that for any $h\in[h_0,h_1]$ we have
\begin{equation}\label{eqnPoisson11}
\left\{\begin{array}{cl}
-\lap\phi_A+A u \cdot \nabla \phi_A  = 0 &\text{in }\Omega \setminus \Omega_h,\\
\phi_A = 0 & \text{ on } \partial \Omega,\\
\ds \int_{\partial \Omega_h} (\nabla \phi_A \cdot \hat n) \, d\sigma =0.
\end{array}\right.
\end{equation}
where the third equation is obtained by integrating the difference of \eqref{aug213} and \eqref{aug214} over $\Omega_h$, and using $u \cdot \hat n = u' \cdot \hat n = 0$ on $\del \Omega_h$. Multiplying (\ref{eqnPoisson11}) by $\phi_A$ and integrating by parts, we obtain:
\[
\int_{\Omega \setminus \Omega_h} \abs{\nabla \phi_A}^2= -\int_{\partial
\Omega_h}\phi_A (\nabla \phi_A \cdot \hat n )\, d \sigma.
\]
Combining the flux condition in \eqref{eqnPoisson11}  with the last
equality, we obtain:
\[
\int_{\Omega \setminus \Omega_h} \abs{\nabla \phi_A}^2 dx = -\int_{\partial
\Omega_h}(\phi_A-\tilde{\phi}_A) \nabla \phi_A \cdot \hat n \, d \sigma,
\]
where
\[
\tilde{\phi}_A\at_{\del \Omega_h} = \frac{1}{\abs{\partial \Omega_h}}\int_{\partial \Omega_h}\phi_A \, d\sigma
\]
is the streamline-averaged $\phi_A$. Integrating this identity for
$h\in[h_0, h_1]$ we obtain:
\begin{align*}
&(h_1-h_0)\int_{\Omega \setminus \Omega_{h_0}} \abs{\nabla \phi_A}^2 \leq
 \int_{h_0}^{h_1} \left( \int_{\Omega \setminus \Omega_h} \abs{\nabla \phi_A}^2
dx \right) \, dh\\
&=-\int_{h_0}^{h_1}\int_{\partial \Omega_h}(\phi_A-\tilde{\phi}_A) (\nabla \phi_A
\cdot \hat n ) \, d\sigma \, dh
\leq \int_{h_0}^{h_1}\int_{\partial \Omega_h}\abs{\phi_A-\tilde{\phi}_A}~\abs{\nabla\phi_A} \, d\sigma \, dh\\
&= \int_{\Omega_{h_0} \setminus \Omega_{h_1}}\!\!\!\!\!
\abs{\phi_A-\tilde{\phi}_A}~\abs{\nabla \phi_A}~\abs{\nabla \psi} \, dx
\leq C\norm{ \nabla \phi_A}_{L^2(\Omega)} \left(\int_{\Omega_{h_0} \setminus \Omega_{h_1}}
\!\!\!\!\!\abs{\phi_A-\tilde{\phi}_A}^2 
dx\right)^{1/2}\!\!\!\!.
\end{align*}
Multiplying \eqref{aug213} by $\tau_A$, \eqref{aug214} by $\tau_A'$, and integrating over $\Omega$, we obtain the uniform bound $\norm{ \nabla \tau_A}_{L^2(\Omega)}, \norm{ \nabla \tau_A'}_{L^2(\Omega)} \leq C$. Hence $\norm{ \nabla \phi_A}_{L^2(\Omega)} \leq C$ and it follows that
\begin{multline}\label{aug218}
\int_{\Omega \setminus \Omega_{h_0}} \abs{\nabla \phi_A}^2 \leq 
C\norm{\phi_A-\tilde{\phi}_A}_{L^2(\Omega_{h_0} \setminus \Omega_{h_1})} \\
\leq C\left(\norm{\tau_A-\tilde{\tau}_A}_{L^2(\Omega_{h_0} \setminus \Omega_{h_1})}+
\norm{\tau_A'-\tilde{\tau}_A'}_{L^2(\Omega_{h_0} \setminus \Omega_{h_1})}\right).
\end{multline}
We claim now that right side of~\eqref{aug218} tends to zero as $A \to \infty$. Indeed, multiplying \eqref{aug213} by $u\cdot\nabla\tau_A$, integrating, using incompressibility of $u$, and the fact that $u\cdot \hat n=0$ on $\partial\Omega$ gives
\begin{eqnarray*}
&&A\int_\Omega(u\cdot\nabla\tau_A)^2dx=\int_\Omega (u\cdot\nabla\tau_A)\lap\tau_A dx
=-\int_\Omega\pdr{\tau_A}{x_j}\pdr{}{x_j}\left(u\cdot\nabla\tau_A\right)dx\\
&&=-\int_\Omega\pdr{\tau_A}{x_j}\pdr{u_m}{x_j} \pdr{\tau_A}{x_m} dx\leq C\int_\Omega
\abs{\nabla\tau_A}^2\leq C.
\end{eqnarray*}
As $\abs{u}$ is strictly positive in $\Omega_{h_0} \setminus \Omega_{h_1}$, it follows that
\begin{equation} \label{2.222}
\norm{\tau_A-\tilde\tau_A}_{L^2(\Omega_{h_0} \setminus \Omega_{h_1})}\to 0
\end{equation}
as $A\to+\infty$. The argument for $\tau_A'$ is identical, completing the proof of \eqref{aug216}.

In order to improve the $\dot H^1(\Omega\setminus\Omega_{h_0})$ bound (\ref{aug216})
to a bound in $L^\infty(\Omega\setminus\Omega_{h_0})$ we simply note that,
given any $\eps>0$ and $A>A_0$,
using (\ref{2.222}) we may find a streamline $\partial\Omega_{h'}$
with $h'<h_0$ but arbitrarily close to $h_0$ so that
\[
\|\tau_A-\tilde\tau_A\|_{L^2(\partial\Omega_{h'})}+\|\tau_A'-\tilde\tau_A'\|_{L^2(\Omega_{h'})}
+\|\nabla\tau_A-\nabla\tau_A'\|_{L^2(\partial\Omega_{h'})}<\eps.
\]
It follows that then 
\[
\|\tau_A-\tilde\tau_A\|_{L^\infty(\Omega_{h'})}+
\|\tau_A'-\tilde\tau_A'\|_{L^\infty(\Omega_{h'})}
<C\eps,
\]
and, in addition, $|\tilde\tau_A-\tilde\tau_A'|_{\partial\Omega_{h'}}<C\eps$ because
of (\ref{aug216}) and since $\tau_A=\tau_A'=0$ on $\partial\Omega$. Finally, since
$\tau_A$ and $\tau_A'$ satisfy the same equation outside of $\Omega\setminus\Omega_{h_0}$,
the maximum principle implies that $|\tau_A-\tau_A'|<C\eps$ in $\Omega\setminus\Omega_{h_0}$.


Now assume that $\psi=\tau^0$ maximizes $\|\tau^u\|_{L^\infty}$ (then $u_0=\nabla^\perp\tau^0$ maximizes $Q_\Omega$) but $u_0$ is not a critical point of $I$. 
Then there exists a $C^4$ stream function $\psi'$ on $\Omega$, equal to $\tau^0$ on $\Omega\setminus\Omega_{h_0}$, such that for $w=\nabla^\perp\psi'$ (when restricted to $\Omega_{h_0}$),
\begin{equation}\label{aug2110}
M=\norm{\tau^0}_{L^\infty}<h_0+Q_{\Omega_{h_0}}(w).
\end{equation}
We can assume that  $\psi'$ has a single critical point in $\Omega_{h_0}$ because so does $\tau^0$ as well as all the perturbations $\psi^\eps$ considered in the proof of Proposition~\ref{ppnCriticalPt}. Moreover, we can assume $\psi'\equiv\tau^0$ on $\Omega\setminus\Omega_{h_1}$ for some $h_1>h_0$ because it is sufficient to consider such perturbations in that proof (see the remark after the proof of Lemma \ref{lmaCriticalPtUnconstrained}).

Then the previous argument shows 
$\norm{\tau^{Aw}-\tau^0}_{L^\infty(\Omega_{h_0-\epsilon}\setminus\Omega_{h_0+\epsilon})}\to 0$ as $A\to+\infty$.
In particular, $\tau^{Aw}>h_0-\delta$ on $\partial\Omega_{h_0}$ for all large $A$, with $\delta=(h_0+Q_{\Omega_{h_0}}(w)-M)/2$. 
This means that $\tau^{Aw}>h_0-\delta+\tau^{Aw}_{\Omega_{h_0}}$ on $\Omega_{h_0}$ by the maximum principle. 
But then
\[
M\ge Q_\Omega(w) \ge h_0-\delta+ Q_{\Omega_{h_0}}(w) =M+\delta >M,
\]
a contradiction. This finishes the proof.
\end{proof}

\subsection{Proof of Theorem~\ref{T.1.1} }  \label{ss2.2}
 
We can assume that $v$ is sufficiently smooth (and approximate general $v$ with smooth ones). Let us denote $\tau=\tau^v$ and $\Omega_h=\Omega_{\tau,h}$. Then by Sard's theorem the set $\mathcal A$ of regular values of $\tau$ has full measure. Thus $\partial\Omega_h$ is a finite union of sufficiently smooth compact manifolds without boundary for each $h\in \mathcal A$ (moreover, $\mathcal A$ is then open because $\tau\in C^2(\Omega)$).

Let $\Omega^*$ and $\tau^*$ be the symmetric rearrangements of $\Omega$ and $\tau$. That is, $\Omega^*$ is the ball with volume $\abs{\Omega}=V$  
 centered at the origin and $\tau^*:\Omega^*\to R_+$ is the non-increasing radial function such that the ball
$\Omega^*_h =  \{ x\in\Omega \,|\, \tau^*(x)>h \}$ satisfies $\abs{\Omega^*_h}=\abs{\Omega_h}$ for each $h\in \mathbb{R}$ (with $\Omega_h$ as above). 

Let now $h\in \mathcal A$. The isoperimetric inequality gives
\begin{equation} \label{1.1}
\abs{\partial\Omega^*_h} \leq \abs{\partial\Omega_h},
\end{equation}
with equality precisely when $\Omega_h$ is a ball. Since $v$ is divergence-free and $\tau$ is constant on $\partial\Omega_h$, we have
\begin{equation} \label{1.2}
\int_{\partial \Omega_h} \abs{\nabla \tau} d\sigma =  -\int_{\partial \Omega_h} \frac{\partial \tau}{\partial \nu} d\sigma= \int_{\Omega_h} (-\lap \tau +v\cdot \nabla \tau)dx = \abs{\Omega_h}=\abs{\Omega^*_h}.
\end{equation}
Finally, the co-area formula yields
\begin{equation} \label{1.3}
-\int_{\partial\Omega_h} \frac{1}{ \abs{\nabla \tau}} \, d\sigma = \frac \partial{\partial h} \abs{\Omega_h} =  \frac \partial{\partial h} \abs{\Omega^*_h} = -\int_{\partial\Omega^*_h} \frac{1}{ \abs{\nabla \tau^*}} \, d\sigma.
\end{equation}
Thus by \eqref{1.1} and the Schwarz inequality,
\[
\int_{\partial \Omega^*_h} \abs{\nabla \tau^*} \, d\sigma \int_{\partial\Omega^*_h} \frac{1}{ \abs{\nabla \tau^*}} \,d\sigma = \abs{\partial\Omega^*_h}^2 \leq \abs{\partial\Omega_h}^2 \leq  \int_{\partial \Omega_h} \abs{\nabla \tau} \,d\sigma \int_{\partial\Omega_h} \frac{1}{ \abs{\nabla \tau}} \,d\sigma.
\]
In view of \eqref{1.2} and \eqref{1.3} we obtain 
\[
\int_{\partial\Omega^*_h}  \abs{\nabla \tau^*} \,d\sigma \leq \int_{\partial\Omega_h}  \abs{\nabla \tau} \, d\sigma = \abs{\Omega^*_h},
\]
with equality precisely when $\Omega_h$ is a ball and $\abs{\nabla \tau}=-\partial \tau/\partial \nu$ is constant on $\partial\Omega_h$. 

So if $\gamma(\abs{x})=\tau^*(x)$ and $\rho=(V/\Gamma_n)^{1/n}$ is the radius of $\Omega^*$, then with $\Sigma_n=n\Gamma_n$ the surface of the unit sphere,
\[
\gamma(\rho)=0 \qquad \text{and} \qquad 0\leq -\gamma'(r)\leq \frac{\Gamma_n r^n}{\Sigma_n r^{n-1}} = \frac r n
\] 
when $\gamma(r)\in \mathcal A$. 
Since $\mathcal A$ has full measure and $\gamma$ is continuous, we have
\begin{equation} \label{1.4}
\gamma(r)\leq \frac{\rho^2-r^2}{2n} = \tilde \gamma(r)
\end{equation}
for all $r\in[0,\rho]$, with $\gamma\equiv \tilde \gamma$ precisely when all $\Omega_h$ are balls and $\tau$ is radial (thus so is $v(x)\cdot x$, hence $v(x)\cdot x\equiv 0$ since $v$ is divergence-free).  Now \eqref{1.4} gives 
\[
\abs{\Omega_h}=\abs{\Omega^*_h}=\abs{\{x\in\Omega^* \,|\, \gamma(\abs{x})>h \}} \leq \abs{\{x\in\Omega^* \,|\, \tilde\gamma(\abs{x})>h \}} \]
and the claim follows. \hfill $\Box$

\section{Properties of the Maximizer}\label{maximizer}
We start by proving \eqref{aug182}.

\begin{proposition}\label{ppnFreidlinProblem}
Let $\psi>0$ be a $C^4$ stream function on a bounded simply connected domain $\Omega\subset\R^2$ with a single 
critical point and let $u = \gradperp \psi$. Then $\tau^{Au} \to \bar \tau^u$ uniformly on $\Omega$, where $\bar\tau^u$ is given by
\begin{equation}\label{eqnFreidlinProblem}
\bar \tau^u(y) = - \int_0^{\psi(y)} \frac {\abs{\Omega_{\psi,h} }} {\int_{\Omega_{\psi,h}} \lap\psi \, dx} \, dh.
\end{equation}
\end{proposition}

\begin{proof}
Assume first that the maximum of $\psi$ is non-degenerate and let $\sup \psi=M>0$. It is then proved in \cite{BKNR}
that, as $A \to \infty$, the functions $\tau^{Au}$ converge uniformly on $\Omega$ to $\bar\tau^u$ with $\bar\tau^u(y)=\bar\tau(\psi(y))$, where $\bar\tau$ solves  the effective problem 
\begin{equation}\label{fr-3}
\left\{
\begin{aligned}
&-\frac{1}{T(h)}\frac{d}{dh}\left(p(h)\frac{d\bar \tau}{dh}\right)=1,\\
&\bar \tau(0)=0 \text{ and $\bar \tau$ is bounded on $(0,M)$}
\end{aligned}
\right.
\end{equation}
on the interval $(0,M)$, with the coefficients
\begin{equation}\label{fr-one-coeff}
T(h)=\int_{\partial \Omega_{\psi,h}}\frac{d\sigma}{\abs{\nabla\psi}},
\quad p(h)=\int_{\partial \Omega_{\psi,h}}{\abs{\nabla\psi}} \, d\sigma.
\end{equation}
By Green's formula and~(\ref{fr-one-coeff}),
\[
p(h)= \int_{\partial \Omega_{\psi,h}}{\abs{\nabla\psi}} \,d\sigma = -\int_{ \partial \Omega_{\psi,h}}{\nabla \psi \cdot \hat n} \, d\sigma
=-\int_{\Omega_{\psi,h} }\lap \psi\, dx,
\]
where we have used $\nabla \psi \cdot \hat n = -\abs{\nabla\psi}$ on $\partial\Omega_{\psi,h}$.
By the co-area formula,
\[
\int_h^{M} T(r) \, dr = \int_{\Omega_{\psi,h} } \, dx,
\]
and thus~\eqref{fr-3} reduces to
\begin{equation}\label{fr-3__}
\bar \tau'(h)=
\frac{\int_h^{M} T(r) \, dr + C}{ p(h)} = 
\frac{\abs{\Omega_{\psi,h}} + C}{-\int_{\Omega_{\psi,h} }\lap \psi \,dx}.
\end{equation}
Non-degeneracy of the maximum of $\psi$ shows that $\tfrac 1{h-M} \int_{\Omega_{\psi,h} }\lap \psi \,dx$ stays bounded away from zero and infinity as $h\uparrow M$. Boundedness of $\bar\tau$ then forces $C=0$,
completing the proof of the non-degenerate case.

If the maximum of $\psi$ is degenerate, we let $\psi_n$ be a $C^4$ stream function with a single non-degenerate critical point which agrees with $\psi$ on $\Omega\setminus \Omega_{\psi,M-(1/n)}$. The proof of Lemma \ref{lemaug212}, with $u=\nabla^\perp\psi$ and $u'=u_n=\nabla^\perp\psi_n$, shows that $\tau^{Au}-\tau^{Au_n}\to 0$ as $A\to\infty$, uniformly on $\Omega_{\psi,M-(2/n)}$. But 
\[
\int_{\Omega_{\psi,h} }\lap \psi \,dx = -\int_{\partial \Omega_{\psi,h}}{\abs{\nabla\psi}} \,d\sigma
\]
shows that $\bar\tau^u$ and $\bar\tau^{u_n}$ coincide on $\Omega_{\psi,M-(1/n)}$, so $\tau^{Au}\to \bar\tau^u$ as $A\to\infty$ uniformly on $\Omega_{\psi,M-(2/n)}$. The result now follows by taking $n\to\infty$ and noticing that for large $n$ and large $A$, the oscillation of $\tau^{Au}$ on $\Omega_{\psi,M-(2/n)}$ has to be small thanks to the small oscillation of $\tau^{Au}$ on $\partial\Omega_{\psi,M-(2/n)}$, small diameter of $\Omega_{\psi,M-(2/n)}$, and the maximum principle.
\end{proof}

We are presently unable to analytically prove existence of solutions to~\eqref{eqnCriticalPt}. The structure of the nonlinear term in~\eqref{eqnCriticalPt} yields itself naturally to
some apriori estimates. These, however, are not strong enough to prove existence, mainly because they do not seem to provide any form of compactness.  
\begin{proposition}\label{ppnApriori}
Let $\phi$ be a $C^4$ solution of~\eqref{eqnCriticalPt} with a single critical point and $\phi = 0$ on $\del\Omega$. Then
\begin{enumerate}
\item $\linfnorm{\phi} \leq \displaystyle\frac{\abs{\Omega}}{4\pi}$.
\item For any Borel function $f$,
$$
\int_{\Omega_{\psi,h}} f(\phi) \, \abs{\lap \phi} \, dx = \int_{\Omega_{\psi,h}} f(\phi) \, dx,
$$
and, in particular, $\displaystyle\int_\Omega \abs{\lap \phi} \, dx= \abs{\Omega}$.
\item If $\tau$ satisfies $-\lap \tau = 1$ in $\Omega$ with $\tau = 0$ on $\del \Omega$, 
then
$$
\int_Omega \abs{\lap \phi - \lap \tau} < \abs{\Omega}.
$$
\end{enumerate}
\end{proposition}

These estimates  do give us some insight as to the nature of classical solutions to~\eqref{eqnCriticalPt}. For instance,
the first two assertions give $L^\infty(\Omega)$ and $H^{1}(\Omega)$ bounds
on $\phi$, while the third is an explicit upper bound 
on the distance between a classical solution of~\eqref{eqnCriticalPt} and the exit time of the Brownian motion from $\Omega$.

\begin{proof}
The second assertion follows by multiplying~\eqref{eqnCriticalPt} by $f(\phi)$ 
and using the co-area formula. As a consequence, for any $h>0$ we have the identity
\[
\abs{\Omega_{\psi,h}}=\int_{\partial\Omega_{\psi,h}}\abs{\nabla\phi} \, d\sigma.
\]
Then (1) follows by
a rearrangement argument as in the proof of Theorem~\ref{T.1.1}.
For the last claim, note that
$$
2(\lap \tau(x) - \lap \phi(x)) = \frac{\abs{\grad \phi(x)}^2}{\ds\int\limits_{\Omega_{\phi,\phi(x)}} 
\abs{\grad \phi} \, d\sigma } \;\int\limits_{\Omega_{\phi,\phi(x)}} 
\frac{1}{\abs{\grad \phi}} \, d\sigma - 1.
$$
By the co-area formula, the integral over $\Omega$ of the first term is exactly $\abs{\Omega}$. Since that term is non-negative, the strict inequality in (3) follows.
\end{proof}

Since an analytical proof of existence of solutions 
for~\eqref{eqnCriticalPt} is at present intangible, we turn our attention to numerics. As boundary integrals are problematic to compute numerically, it is more convenient to work with
the equation
\begin{equation}\label{eqnCriticalPtLn}
-2\lap \phi(x) = 1 + \grad \phi(x) \cdot \grad \ln \abs{\bigl.\Omega_{\phi(x)}},
\end{equation}
which is equivalent to~\eqref{eqnCriticalPt}.  Surprisingly, an iteration scheme of the form
\begin{gather*}
-\lap \phi_0 = 1, \quad \phi\at_{\del\Omega} = 0\\
-2 \lap \phi_{n+1}(x) = 1 + \grad \phi_n(x) \cdot \grad 
\ln \abs{ \bigl. \{\phi_n \geq \phi_n(x) \} }, \quad \phi_{n+1} \at_{\del \Omega} = 0
\end{gather*}
does not always converge. For certain domains, it turns out that numerically 
$\linfnorm{\phi_n} \to \infty$ as $n \to \infty$, which is clearly not representative 
of the solution of~\eqref{eqnCriticalPt} as it violates the last assertion in 
Proposition~\ref{ppnApriori}.

\begin{figure}[thb]\label{fgrMaximisersAndExitTimes}
  \subfigure[Maximiser $\psi$]{
    \includegraphics[width=5cm]{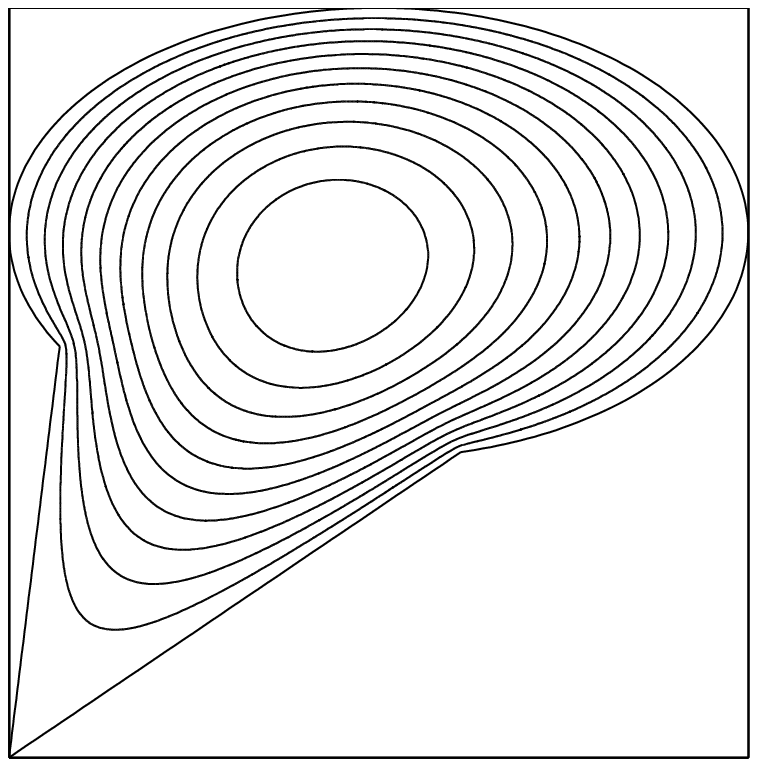}
  }
  \quad
  \subfigure[Expected exit time $\tau_0$]{
    \includegraphics[width=5cm]{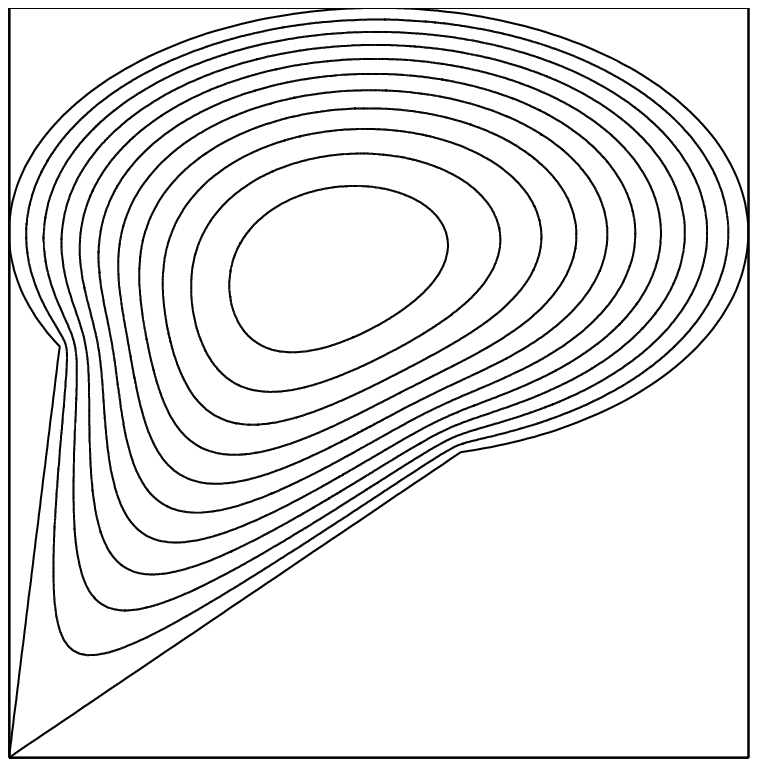}
  }\\
  \subfigure[Maximiser $\psi$]{
    \includegraphics[width=5cm]{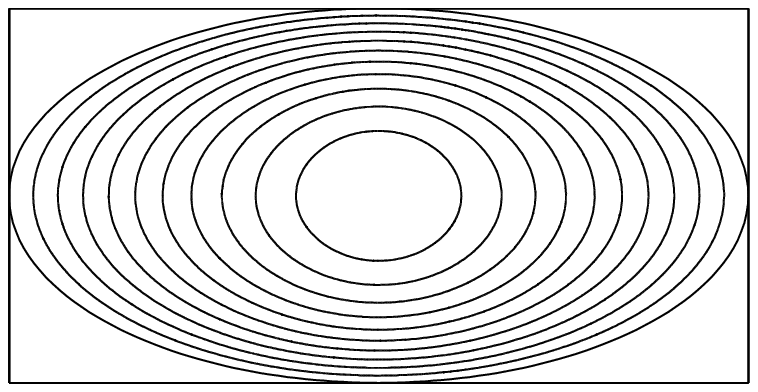}
  }
  \quad
  \subfigure[Expected exit time $\tau_0$]{
    \includegraphics[width=5cm]{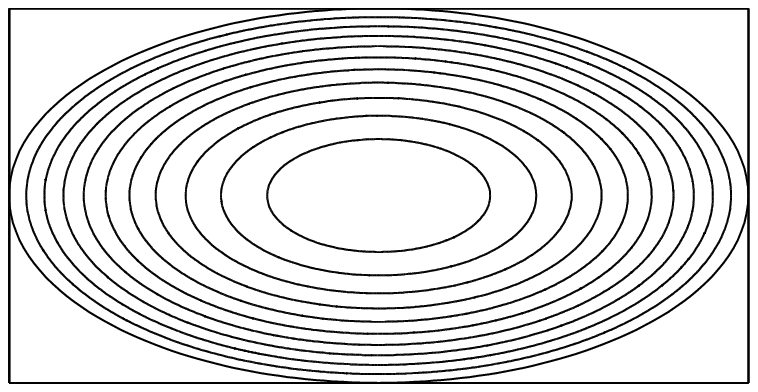}
  }\\
\caption{Maximisers and expected exit times from two different domains.}
\end{figure}

It turns out that an iteration scheme of the form
\begin{gather}
-\lap \phi_0 = 1, \quad \phi\at_{\del\Omega} = 0\\
-2 \lap \phi_{n+\frac{1}{2}}(x) = 
1 + \grad \phi_n(x) \cdot \grad 
\ln \abs{ \bigl. \{\phi_n \geq \phi_n(x) \} }, \quad \phi_{n+\frac{1}{2}}\at_{\del\Omega}=0\\
\label{eqnIterScheme3} \phi_{n+1}\at_{\phi_{n+\frac{1}{2}} = h_0} = 
- \int_0^{h_0} \frac{\abs{\bigl\{\phi_{n+\frac{1}{2}} \geq h\bigr\} }}{\int_{\{\phi_{n+\frac{1}{2}} = h\}} \frac{\del \psi}{\del \hat n} \, d\sigma} \, dh.
\end{gather}
does converge rapidly to a numerical solution of~\eqref{eqnCriticalPtLn}. In fact, \eqref{eqnIterScheme3} can be replaced by
\begin{equation*}\tag{$\text{\ref{eqnIterScheme3}}$'}
-\lap \phi_{n+1} + A \gradperp \phi_{n+\frac{1}{2}} \cdot \grad \phi_{n+1} = 1
\end{equation*}
for some large, fixed $A$, which produces better numerical results. Figure~\ref{fgrMaximisersAndExitTimes} shows contour plots of the solution to~\eqref{eqnCriticalPt} in two different domains. For comparison, the expected exit time from the domain $\tau_0$ is shown alongside each plot of $\phi$.
 
We are unable to prove convergence of these numerical schemes, just as we can not establish existence of solutions of~\eqref{eqnCriticalPt}. However, one immediate observation from Figure~\ref{fgrMaximisersAndExitTimes} is that the level sets of $\phi$ become circular near the maximum. Indeed, for any classical solution of~\eqref{eqnCriticalPt}, this must be the case.
\begin{proposition}
Let $\phi$ be a smooth solution of~\eqref{eqnCriticalPt} and assume that $\phi$ attains a 
local maximum at $(0,0)$, then $\del_{xx} \phi(0) = \del_{yy} \phi(0)$.
\end{proposition}
\begin{proof}
We will show if $\phi$ is any smooth function which attains a maximum at $0$, then the last term in~\eqref{eqnCriticalPtLn} is continuous near $0$ if and only if $\del_{xx} \phi(0) 
= \del_{yy} \phi(0)$. This immediately implies the proposition.

Assume first that the Hessian of $\phi$ at $0$ is not degenerate (in this case 
assuming $\phi \in C^3$ will be enough for the proof). We rotate our coordinate frame and assume without loss of generality that
$$
\phi(x,y) = M - \frac{x^2}{a^2} - \frac{y^2}{b^2} + c_3(x,y)
$$
where $c_3(x,y)$ is some function involving only third order or higher terms. Now for any $\epsilon > 0$, define $f(\epsilon)$ by
\begin{equation}
f(\epsilon) = \abs{ \Omega_{\phi, M - \epsilon} } = \abs{\left\{ \frac{x^2}{a^2} + \frac{y^2}{b^2} \leq \epsilon + c_3(x,y) \right\}}
\label{eqnFofEpsilon}
= \pi a b \epsilon + O(\epsilon^2),
\end{equation}
whence
\begin{equation}
\label{eqnLnFPrime}
 \frac{d}{d\epsilon} \ln (f(\epsilon)) = \frac{1}{\epsilon} + O(1).
\end{equation}
Thus
\begin{align*}
\grad \phi(x,y) \cdot \grad \ln \abs{ \Omega_{\phi, \phi(x,y)}} 
&= \abs{\grad \phi(x,y)}^2 \left( \frac{1}{M - \phi(x,y)} + O(1) \right)\\
    &= \frac{\abs{\grad \phi(x,y)}^2}{ M - \phi(x,y)} + O(1)\abs{\grad \phi(x,y)}^2
\end{align*}
The second term on the right is certainly continuous at $0$, since $\grad \phi(0) = 0$. The first term is continuous at $0$ if and only if $a = b$.

In the case that the Hessian of $\phi$ is degenerate at $0$, 
the above proof works with minor modifications. Using a higher order Taylor approximation of $\phi$, the right hand side of \eqref{eqnFofEpsilon} becomes $c_1 \epsilon^{c_2}$ for two constants $c_1 > 0$, and $c_2 = \frac{1}{2} + \frac{1}{n}$ with $n \geq 3$. 
Now replacing ${1}/{\epsilon}$ with 
${c_2}/{\epsilon}$ in~\eqref{eqnLnFPrime}, the remainder of the proof is unchanged.
\end{proof}

\section{Proof of Proposition~\ref{ppnCriticalPt}}\label{ptproof}

First, we obtain an expression for $V(\psi,v)$. Let $\Omega_h=\Omega_{\psi,h}$ and $\Omega_h^\epsilon = \Omega_{\psi^\epsilon, h } = X_\epsilon \inv(\Omega_h)$.

\begin{lemma}\label{lmaVariationTau}
Let $M=\sup\psi$, then the variation~\eqref{eqnVariationTau} is
\begin{multline}\label{eqnDDEpsilon}
V(\psi, v ) =
\int_0^M \frac{1}{\left(\int_{\del\Omega_h} \frac{\del \psi}{\del \hat n} \, d\sigma\right)^2}
\Biggl[ \abs{\Omega_h}\int_{\del \Omega_h} \left( \frac{\del}{\del \hat n} \left(v \cdot \grad \psi \right) - \lap \psi\,  v \cdot \hat n \right)\, d\sigma +\\
+ \left(\int_{\del \Omega_h} \frac{\del \psi}{\del \hat n} \, d\sigma \right) \left( \int_{\del \Omega_h} v \cdot \hat n \, d\sigma\right)\Biggr] \, dh.
\end{multline}
\end{lemma}

\begin{proof}
This follows from Proposition~\ref{ppnFreidlinProblem}. Note that
\begin{multline}\label{eqnDDEpsilonAbsOmegah}
\left. \frac{d}{d\epsilon} \right\rvert_{\epsilon = 0} \abs{\Omega_h^\epsilon} = \left. \frac{d}{d\epsilon} \right\rvert_{\epsilon = 0} \int_{\Omega_h^\epsilon} dx   = \left. \frac{d}{d\epsilon} \right\rvert_{\epsilon = 0} \int_{\Omega_h} \abs{\det \grad X_\epsilon \inv} \, dx   \\
= -\int_{\Omega_h} \divergence v dx = -\int_{\del \Omega_h} v \cdot \hat n \, d\sigma,
\end{multline}
and
\begin{align}
\nonumber \left.\frac{d}{d\epsilon}\right\rvert_{\epsilon = 0} \int_{\del \Omega_h^\epsilon} \frac{\del \psi^\epsilon}{\del \hat n} \, d\sigma &= \left.\frac{d}{d\epsilon}\right\rvert_{\epsilon = 0} \int_{\Omega_h^\epsilon} \lap \psi^\epsilon \, dx \\
\nonumber &= \left.\frac{d}{d\epsilon}\right\rvert_{\epsilon = 0} \int_{\Omega_h} \left(\lap \psi^\epsilon\right) \circ X_\epsilon\inv \, \abs{\det \grad X_\epsilon\inv} \, dx \\
\nonumber &= \int_{\Omega_h} \Bigl[ \lap \left( v \cdot \grad \psi \right) - v \cdot \grad \lap \psi - \left(\divergence v\right) \lap \psi \Bigr] \, dx \\
\label{eqnDDEpsilonDPhiDn} &= \int_{\del \Omega_h}  \frac{\del}{\del \hat n} \left(v \cdot \grad \psi\right) \, d\sigma - \int_{\del \Omega_h} (\lap \psi) v \cdot \hat n \, d\sigma
\end{align}
Thus, using~\eqref{ppnFreidlinProblem}, and equations~\eqref{eqnDDEpsilonAbsOmegah}--\eqref{eqnDDEpsilonDPhiDn} we are done.
\end{proof}

Before proving Proposition~\ref{ppnCriticalPt}, we require a lemma.

\begin{lemma}\label{lmaCriticalPtUnconstrained}
A $C^4$ stream function $\psi$ (with a single critical point) is a critical point of the functional $I$ if and only if it solves 
\begin{equation}\label{eqnCriticalPtVarphi}
\grad F_\psi \cdot \grad \psi + 2 F_\psi \lap \psi - G_\psi = 0,
\end{equation}
where $F_\psi$ and $G_\psi$ are defined by
\begin{eqnarray*}
G_\psi(x) = \left( \int\limits_{\Omega_{\psi(x)}} \frac{\del \psi}{\del \hat n} \, d\sigma \right)\inv
\quad\text{and}\quad
F_\psi(x) = \abs{\bigl.\Omega_{\psi(x)}} G_\psi(x)^2.
\end{eqnarray*}
\end{lemma}
\begin{proof}
With $F_\psi,G_\psi$ as above, equation~\eqref{eqnDDEpsilon} reduces to
\begin{equation}\label{eqnVFG}
V(\psi, v) = \int_0^M F_\psi \int\limits_{\partial\Omega_{h}} 
\left[ \frac{\del}{\del \hat n} \left( v \cdot \grad \psi \right) - \lap \psi \, v \cdot \hat n \right] \, d\sigma \, dh 
+ \int_{0}^M G_\psi \int\limits_{\partial\Omega_{h}} 
v \cdot \hat n \, d\sigma \, dh.
\end{equation}
By the co-area formula, we have, first,
\begin{equation}\label{aug242}
\int_{0}^M G_\psi \int\limits_{\partial\Omega_{h}} 
v \cdot \hat n \, d\sigma \, dh 
= \int_0^M \int\limits_{\partial\Omega_{h}} 
(-G_\psi\, v \cdot \grad \psi) \, \frac{d\sigma}{\abs{\grad \psi}} \, dh
= -\int_\Omega G_\psi \, v\cdot \grad \psi \, dx ,
\end{equation}
second,
\begin{equation}\label{eqnVFGline1term2}
\int_0^M F_\psi  \int\limits_{\partial\Omega_{h}} 
\left[- \lap \psi \, v \cdot \hat n \right] \, d\sigma \, dh 
= \int_\Omega F_\psi \, \lap \psi \, (v \cdot \grad \psi) \, dx ,
\end{equation}
and, finally,
\begin{align}
\nonumber \int_0^M F_\psi  \int\limits_{\partial\Omega_{h}} \frac{\del}{\del \hat n} \left( v \cdot \grad \psi \right) \, d\sigma \, dh &= -\int_\Omega F_\psi \, \grad( v \cdot \grad \psi ) \cdot \grad \psi \, dx \\
\label{eqnVFGline1term1} &= \int_\Omega (v\cdot\grad \psi) \left( \grad F_\psi \cdot \grad \psi + F_\psi \lap \psi \right)\, dx .
\end{align}
In the last equality we used the identity
$$
(v\cdot\grad \psi) \left( \grad F_\psi \cdot \grad \psi + F _\psi \lap \psi \right)
+F_\psi \, \grad( v \cdot \grad \psi ) \cdot \grad \psi = \divergence \left[ F _\psi \, (v \cdot \grad \psi) \, \grad \psi \right]
$$
and the fact that $F _\psi \, (v \cdot \grad \psi) \, \grad \psi = 0$ on $\del \Omega$.

Using~\eqref{aug242}--\eqref{eqnVFGline1term1}, expression~\eqref{eqnVFG} becomes
$$
V(\psi, v) = \int_\Omega (v\cdot \grad \psi) \left[\grad F_\psi \cdot \grad \psi + 2 F_\psi \lap \psi - G_\psi \right] \, dx .
$$
Thus $V(\psi, v) = 0$ for all $C^4$ functions $v$ which vanish on $\del \Omega$ if and only if equation~\eqref{eqnCriticalPtVarphi} holds.
\end{proof}

Notice that the same conclusion is obtained if we ask $V(\psi, v)=0$ only for all $v$ compactly supported inside $\Omega$.

\begin{proof}[Proof of Proposition~\ref{ppnCriticalPt}]
Note that the variation $V(\psi, v)$ in~\eqref{eqnVariationTau} depends only on the geometry of the level sets of the stream function 
$\psi$ and thus is invariant under reparametrizations. Thus if $\psi$ is a solution of~\eqref{eqnCriticalPtVarphi}, then for any monotone function $f$, $f\circ\psi$ is also a solution of~\eqref{eqnCriticalPtVarphi}. Note that $\phi = \bar \tau^{\gradperp \psi}$, the solution of the Freidlin problem~\eqref{aug182} with stream function $\psi$,  is only a reparametrization of the level sets of $\psi$. Thus to prove Proposition~\ref{ppnCriticalPt} we only need to show that if $\psi$ solves~\eqref{eqnCriticalPtVarphi}, then $\phi$ solves~\eqref{eqnCriticalPt}.

Note that since $\phi$ solves the Freidlin problem~\eqref{aug182}, we have
\begin{equation}\label{eqnFriedlinConstraint}
\abs{\Omega_{\phi,h}} = \int\limits_{\partial\Omega_{\phi,h}} \abs{\grad \phi} \, d\sigma,
\end{equation}
and so $F_\phi = -G_\phi$. We also have
\begin{align*}
\grad F_\phi(x) &= -\grad G_\phi(x) = G_\phi(x)^2\, \grad\left( 
\int\limits_{\partial\Omega_{\phi,\phi(x)}} -\abs{\grad \phi} \, d\sigma\right)\\
&= -G_\phi(x)^2\, \grad \abs{\bigl.\Omega_{\phi,\phi(x)}} = -G_\phi(x)^2 \, \grad \left( \int\limits_{\phi(x)}^M \int\limits_{\partial\Omega_{\phi,h}} \frac{1}{\abs{\grad \phi}} \, d\sigma \, dh \right)\\
&= -G_\phi(x)^2 \, \grad \phi (x) \, \int\limits_{\partial\Omega_{\phi,\phi(x)}} \frac{1}{\abs{\grad \phi}} \, d\sigma,
\end{align*}
and using this in~\eqref{eqnCriticalPtVarphi} immediately yields~\eqref{eqnCriticalPt}.
\end{proof}

We remark that any solution to~\eqref{eqnCriticalPt} is automatically a solution to the Freidlin problem with itself as stream function (i.e.\ satisfies~\eqref{eqnFriedlinConstraint}). Indeed, integrating~\eqref{eqnCriticalPt} over $\Omega_{\phi,h_0}$ and using the co-area formula gives
$$
-2\int\limits_{\Omega_{\phi, h_0}} \lap \phi = \abs{\bigl.\Omega_{\phi,h_0}} + 
\int_{h_0}^M \int\limits_{\partial\Omega_{\phi,h}} 
\abs{\grad \phi}^2 \frac{d\sigma}{\abs{\grad \phi}} 
\int\limits_{\partial\Omega_{\phi,h}} \frac{d\sigma}{\abs{\grad \phi}}
\left(\,\int\limits_{\partial\smash{\Omega_{\phi,h}}} \abs{\grad \phi} \, d\sigma\right)^{-1} dh,
$$
and hence
$$
2 \int\limits_{\partial\Omega_{\phi,h_0}} \abs{\grad \phi} = \abs{\bigl.\Omega_{\phi,h_0}} 
+ \int_{h_0}^M \int\limits_{\partial\Omega_{\phi,h}} 
\frac{1}{\abs{\grad \phi}} \, d\sigma \, dh =2 \abs{\bigl.\Omega_{\phi, h_0}},
$$
showing~\eqref{eqnFriedlinConstraint} is satisfied.

\end{document}